\documentclass[fleqn,10pt]{SelfArx} 

\definecolor{color1}{RGB}{0,0,90} 
\definecolor{color2}{RGB}{0,20,20} 

\newlength{\tocsep}
\usepackage{lipsum} 
\usepackage{multicol}
\usepackage{algorithmicx}
\usepackage{algorithm}
\usepackage{blindtext}
\usepackage{algpseudocode}
\usepackage{amsmath,amssymb}
\usepackage{varwidth}
\usepackage{lineno}
\JournalInfo{2014} 
\Archive{\textit{Preprint}} 

\PaperTitle{The Plant Propagation Algorithm: Modifications and Implementation} 

\Authors{Muhammad Sulaiman\textsuperscript{1,2}*, Abdellah Salhi\textsuperscript{1}, Eric S Fraga\textsuperscript{3}} 
\affiliation{\textsuperscript{1}\textit{Department of Mathematical Sciences, University of Essex, Colchester CO4 3SQ, UK}} 
\affiliation{\textsuperscript{2}\textit{Department of Mathematics, Abdul Wali Khan University Mardan, KPK, Pakistan}} 
\affiliation{\textsuperscript{3}\textit{Centre for Process Systems Engineering, Department of Chemical Engineering, University College London (UCL)}} 
\affiliation{*\textbf{Corresponding author}: sulaiman513@yahoo.co.uk} 

\Keywords{Optimisation, Heuristics, Nature-Inspired Algorithms, Artificial Bee Colony Algorithm, Strawberry/ Plant Propagation Algorithm} 
\newcommand{\keywordname}{Keywords} 

\Abstract{The Plant Propagation Algorithm, epitomised by the Strawberry Algorithm, has been previously successfully tested on low dimensional continuous optimisation problems. It is a neighbourhood search algorithm.  In this paper, we introduce, robust and efficient version of the algorithm and explain how it can be implemented to compete with one of the best available alternatives, namely the Artificial Bee Colony algorithm and we present an improved and more effective variant on standard continuous optimisation test problem instances in high dimensions. Computational and comparative results are included.}

\begin{document}

\flushbottom 

\maketitle 


\thispagestyle{empty} 


\section*{Introduction} 

\addcontentsline{toc}{section}{\hspace*{-\tocsep}Introduction}

Although there are already many good algorithms and heuristics for optimisation and search problems \cite{talbi2009metaheuristics}, the growing complexity of these problems in practice and the frequency with which they occur mean that new and more effective algorithms have to be developed. Note that frequently occuring problems may justify introducing new algorithms with only slight improvements.  But, designing new algorithms which are easier to implement and require fewer arbitrarily set parameters, for instance, is a worthwhile quest in itself. So, there are many reasons for trying to invent new algorithms.

As it happens, a lot of attempts at creating new algorithms look to Nature for inspiration. It seems that natural phenomena such as the survival of living entities and the success of some species in a given environment, rely on optimisation and search to overcome the constraints that these environments impose on them.  The survival of species often depends on their ability to adapt, to find food quickly, avoid predation and give their off-spring the best chance to survive and thrive. These are typically optimisation/search problems.

The Plant Propagation Algorithm (PPA) is a Nature-inspired algorithm, \cite{yang2011nature, brownlee2011clever}, for optimisation and search. It emulates the way plants, in particular the strawberry plant, propagate. Basic PPA has been described and tested on single objective as well as multi-objective continuous optimisation problems \cite{Salhi2010PPA,sulaiman2014plant}.  Although the problems considered were of low dimensions, it was established that the algorithm has merits and deserves further investigation and testing on higher dimension instances.  The attraction of the algorithm is its simplicity and the relatively small number of parameters requiring arbitrary setting. This paper addresses the issue of testing PPA on larger, higher dimension problems, and compares the algorithm to two other Nature-inspired ones specifically the Artificial Bee Colony (ABC) algorithm \cite{karaboga2005idea} and the Modified Artificial Bee Colony (MABC) algorithm \cite{akay2010modified,gao2011modified}.  The PPA is presented both in its original and improved forms. Extensive comparative results on high-dimensional test instances are reported and discussed. The paper ends with a conclusion and further issues for consideration.


\section{The Artificial Bee Colony Algorithm}

The Artificial Bee Colony algorithm proposed in \cite{karaboga2005idea}, simulates the foraging behaviour of bees living in a colony. Three groups of bees participate in the foraging process: worker bees, onlooker bees and scout bees. The majority of the population is composed of worker bees and onlooker bees. The scouts are recruited from worker bees. Algorithm 1 below describes the ABC algorithm.
\begin{algorithm}\label{ABC}
\caption{\textbf{Outline of the ABC algorithm}, \cite{karaboga2005idea}}
\label{alg:ABC}
\begin{algorithmic}[1]
\State  Initialisation: Generate food sources for all worker bees;
\While{stopping criterion not satisfied}
\State   \begin{varwidth}[t]{\linewidth} Each worker bee goes to its assigned food source, \\finds   a neighbouring source, evaluates it and \\displays   the information to onlooker bees in the hive, through a dance; \end{varwidth}
\State  \begin{varwidth}[t]{\linewidth} Each onlooker bee, after watching worker bees\\ dancing, selects one of the sources and goes to it. \\It then finds a neighbouring source and evaluates it. \end{varwidth}
\State   \begin{varwidth}[t]{\linewidth}Food sources that are not good are abandoned and\\ replaced by new ones discovered by worker bees that have become scouts.\end{varwidth}\vspace{2mm}
\State    Record the best food source found so far.
\EndWhile
\end{algorithmic}
\end{algorithm}
\subsection{The Modified Artificial Bee Colony Algorithm}

The MABC algorithm has been suggested in \cite{gao2011modified}. The new method improves upon the exploitation aspect. MABC uses Differential Evolution \cite{rahnamayan2008opposition} in step 4 of the ABC algorithm, and removes step 6 (or the Scouts phase). Although MABC has a good performance on continuous unconstrained optimisation, it is a hybrid algorithm of ABC. It is, at the moment, the algorithm to beat among a large selection of benchmark functions \cite{gao2011modified}.

\section{The Strawberry Algorithm}

PPA as the Strawberry algorithm, is a neighbourhood search algorithm. However, it can be seen as a multi-path following algorithm unlike Simulated Annealing (SA) \cite{aarts1997simulated, salhi2000experiences,sulaiman2014plant}, for instance, which is a single path following algorithm.

Exploration and exploitation are the two main properties global optimisation algorithms ought to have \cite{yang2011nature,talbi2009metaheuristics,brownlee2011clever,sun2013intelligent}. Exploration refers to the property of covering the whole search space, while exploitation refers to the property of searching for local optima, near good solutions. Effective global optimisation methods exhibit both properties.

Consider what a strawberry plant, and possibly any plant which propagates through runners, will do to maximise its chances of survival. If it is in a good spot in the ground, with enough water, nutrients and light, then it is reasonable to assume that there is no pressure on it to leave that spot to guarantee its survival. So, it will send many short runners that will give new strawberry plants and occupy the neighbourhood as best they can. If, on the other hand, the mother plant is in a spot that is poor in water, nutrients, light or any element necessary for a plant to survive, then it will try to find a better spot for its off-spring. Therefore, it will send a few runners further afield to explore distant neighbourhoods. One can also assume that it will send only a few, since sending a long runner is a big investment for a plant which is in a poor spot.
We may further assume that the quality of the spot (abundance of nutrients, water and light) is reflected in the
growth of the plant. With this in mind, and the following notation, PPA can be described as follows.

A plant $p_i$ is in spot $X_i$ in dimension $n$. This means $X_i=[x_{i,j}], \mbox{ for } j=1,...,n$.
In PPA, exploitation is implemented through sending of many short runners by plants in good spots. Exploration is implemented by sending few long runners by plants in poor spots; the long runners allow distant neighbourhoods to be explored.

The parameters used in PPA are the population size $NP$ which is the number of strawberry plants to start with, the maximum number of generations $g_{max}$, and the maximum number of possible runners $n_{max}$ per plant. $g_{max}$ is effectively the stopping criterion in this initial version of PPA. The algorithm uses the objective function value at different plant positions $X_i, i=1,...,NP$ to rank them as would a fitness function in genetic algorithms, \cite{holland1975adaption}.

Let $N_i\in(0,1)$ be the normalised objective function value for $X_i$.  The number of plant runners, $n_r^i$ for this solution is given by
\begin{equation} \label{eq:numberofrunners}
  n_r^{i}= \lceil n_{max} N_{i}\alpha_i \rceil
\end{equation}
where $\alpha_i \in (0,1)$ is a randomly generated number.
Every solution generates at least one runner.  Each runner generated has a distance, $dx_j^i \in [-1,1]^n$, calculated by
\begin{equation} \label{eq:distanceofrunners}
 dx^i_j=2(1-N_{i})(r-0.5), \mbox{ for } j=1,...,n,
\end{equation}
where $r \in [0,1]$ is also randomly generated.
We note that the number of runners is proportional to the fitness, i.e. the value of the objective function, whereas the distance is inversely proportional to it.

Having calculated $dx^i$, the new point to explore, $Y_i=[y_{i,j}], \mbox{ for } j=1,...,n$, is given by
\begin{equation}
y_{i,j}=x_{i,j}+(b_{j}-a_{j})dx^i_j, \mbox{ for } j=1,\ldots,n
\end{equation}
where $a_j$ and $b_j$ are the lower and upper bounds of the search domain respectively.  If the bounds of the search domain are violated, the point is adjusted to be within the domain.

After all individuals/ plants in the population have sent out their allocated runners, new plants are evaluated and the whole increased population is sorted. To keep the population constant, individuals with lower fitness are eliminated.

\begin{algorithm}\label{PPA}
\caption{\textbf{: Pseudo-code of PPA, \cite{Salhi2010PPA}}}
\label{alg:PPA}
\begin{algorithmic}[1]
\State  Generate a population $P=\{X_i,i=1,...,NP\}$;
\State  $g\leftarrow1$;
   \For {$g=1$ :$g_{max}$}
       \State Compute $N_{i}=f(X_i),\forall \mbox{ } X_{i} \in P$;
       \State \begin{varwidth}[t]{\linewidth} Sort $P$ in ascending order of fitness values $N$ (for minimization); \end{varwidth}
       \State Create new population $\Phi$;
    \For {each $X_i,i=1,...,NP$}

        \State \begin{varwidth}[t]{\linewidth} $r_{i} \leftarrow$ set of runners where both the size of the set \par  \hspace{1mm} \hskip \algorithmicindent   and the distance for each runner (individually) \par  \hspace{2mm}\hskip \algorithmicindent  are  proportional to the fitness values $N_{i}$;\end{varwidth}

     \State \begin{varwidth}[t]{\linewidth} $\Phi \leftarrow$ $\Phi \cup r_{i}$ \{append to population; death occurs \par \hspace{13mm} \hskip \algorithmicindent  by omission\};\end{varwidth}
    \EndFor
    \State $P\leftarrow  \Phi$ \{new population\};
    \EndFor
\State \emph{\textbf{Return: }}$P$, the population of solutions.
   \end{algorithmic}
\end{algorithm}

\section{Modified Version of PPA}

The modifications with respect to the previous version of PPA as in Algorithm 2, \cite{Salhi2010PPA,sulaiman2014plant}, concern the strategy for generating runners, whether they should be short or long, and whether the new points resulting from these runners are retained or not. Also, while in Algorithm~\ref{alg:PPA} the number of runners is between 1 and $n_{max}$, in the new version it is a fixed number, although all generated runners may not lead to points that will make it into the new population. This is because, after sorting, their rank may be above $NP$, the population size.
\subsection{An Alternative Implementation of the Propagation Phase}
The population is initialized randomly by
\begin{equation}x_{i,j}=a_j+(b_j-a_j)\alpha_j,\mbox{ } j=1,...,n\end{equation}
where $\alpha_j \in(0,1)$ is randomly generated for each $j$.
After the population is initialized, MPPA proceeds to generate for every member in the population a number $n_r$ of runners; $n_r$ is a fixed constant. These runners lead to new solutions as per the Equations 5-7, on the off chance that the limits of the search area are maltreated, the  coordinates are conformed respectively to be inside the search space.
\begin{equation}\label{eq:perturbation1}
  y_{i,j}=x_{i,j}+\beta_j x_{i,j}, \qquad j=1,\ldots,n
\end{equation}
where $\beta_j \in [-1,1]$ is a randomly generated number for each $j$. The term $\beta_j x_{i,j}$ is the length with respect to the $j^{th}$ coordinate of the runner, and $y_{i,j}\in[a_j,b_j]$. If the bounds of the search domain are violated, the point is adjusted to be within the domain. The generated individual $Y$ is evaluated according to the objective function and is stored in $\Phi$. Equation (5) helps in exploring the neighbourhood of $x_{i,j}$. As the search becomes refined, that is the algorithm is in exploitation mode then the coordinates produced by Equation (5) are smaller and smaller. This is represented in Figure (1), where the horizontal axis shows the total number of perturbations produced during 30 independent experiments. Note that in beginning of each experiment the step size is larger and as the search is refined the step size decreases gradually, in the later case the algorithm is in exploitation mode. In Algorithm (3), If this newly created solution, by Equation (5), is not improving the objective function, then another individual is created with a runner based on Equation (6). The number of runners of a certain length generated by Equation (6), when solving $f_2$, is shown in Figure (3), this shows the frequency of exploration around the optimum solution.
\begin{equation}\label{eq:perturbation2}
y_{i,j}=x_{i,j}+\beta_j b_j, \qquad j=1,\ldots,n
\end{equation}

\begin{figure*}
\begin{center}
\hspace{-1mm}
\includegraphics{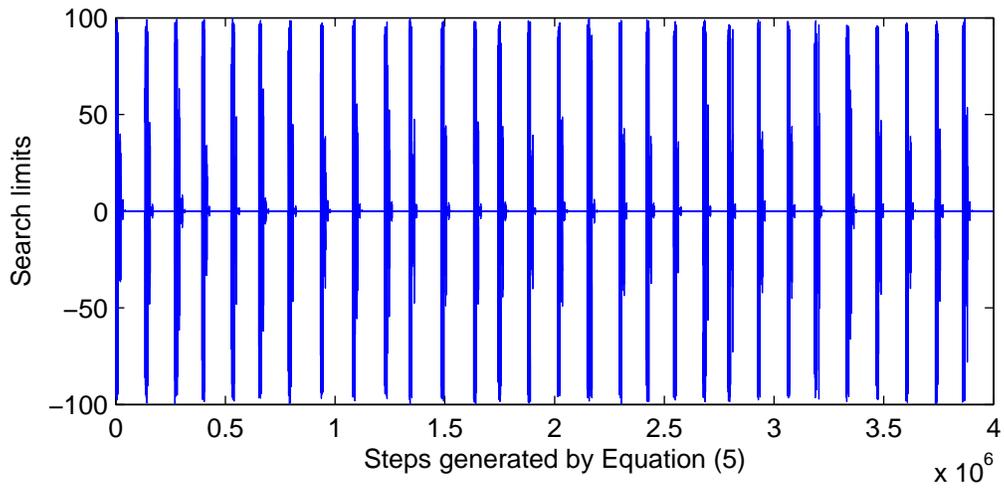}
\caption{Effects of Equation (5) in solving $f_2$}
\end{center}
\end{figure*}

\begin{figure*}
\hspace{40mm}
\includegraphics{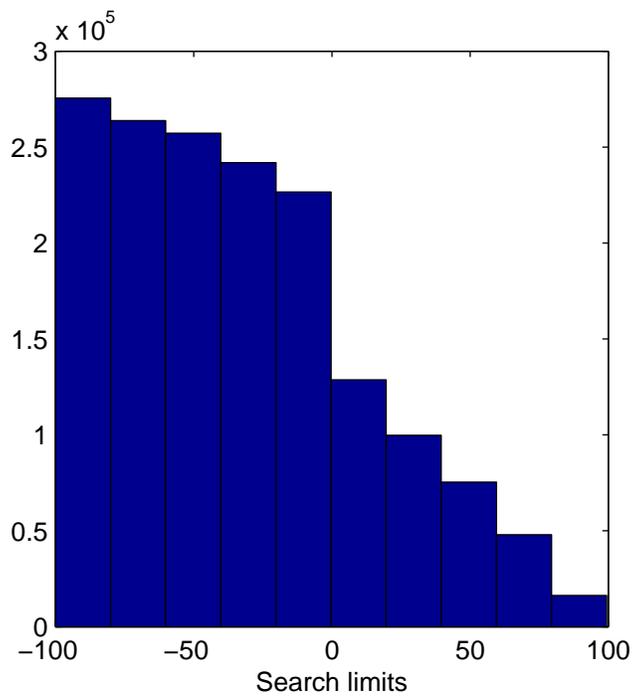}
\caption{Effects of Equation (7) in solving $f_2$}
\end{figure*}

\noindent where $b_j$ is the $j^{th}$ upper bound and here again $y_{i,j}\in[a_j,b_j]$. This can be considered as a solution at the end of a long runner. Again, if the generated individual does not improve the objective value, another runner is created by Equation (7).
\begin{equation}\label{eq:perturbation3}
  y_{i,j}=x_{i,j}+\beta_j a_j, \qquad j=1,\ldots,n
\end{equation}
where $a_j$ is the $j^{th}$ lower bound and $y_{i,j}\in[a_j,b_j]$. This can be considered as a solution at the end of a long runner. The number of runners of a certain length generated by Equation (7), when solving $f_2$, is shown in Figure (2), this shows the frequency of exploration around the optimum solution.

These Equations (5-7), are implemented by Algorithm (3) turn by turn if any of the earlier search equation does not improve the current solution. This improve the balance between exploration and exploitation of the search space as depicted in Figures (1-3). Note that the above equations may lead to infeasibility. In these situations, the offending entry is set by default to the boundary, lower or upper as per the concerned equation.

\begin{figure*}
\hspace{35mm}
\includegraphics{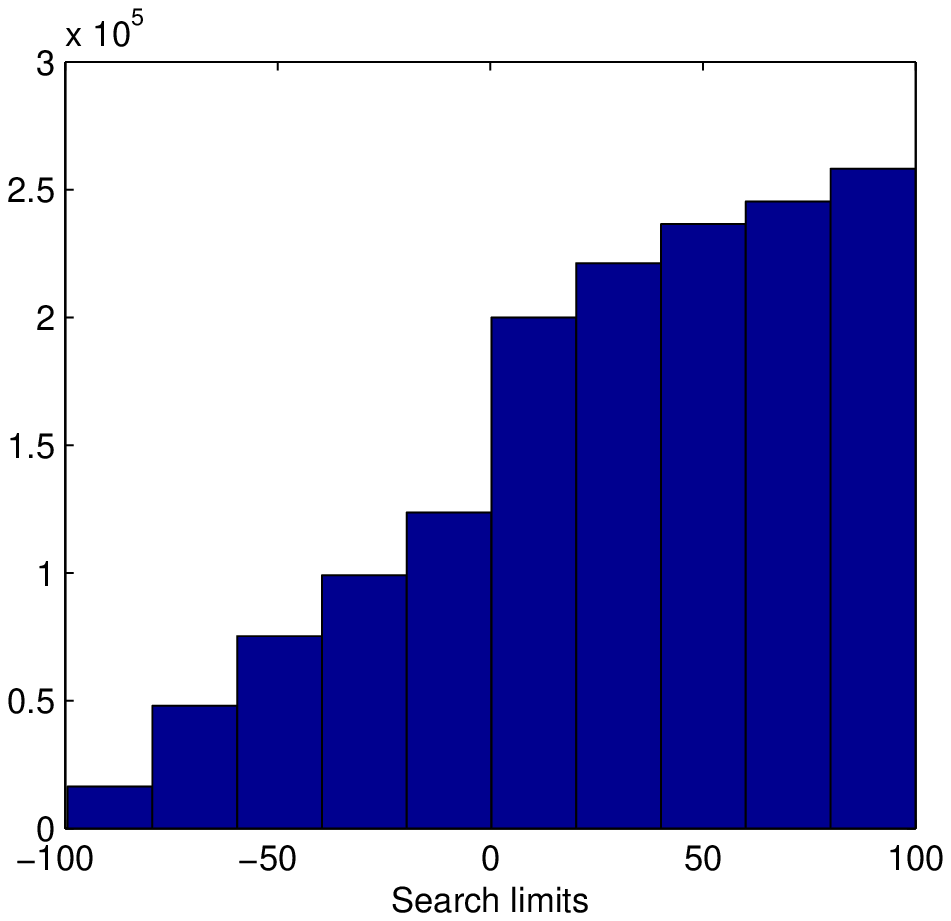}
\caption{Effects of Equation (6) in solving $f_2$}
\end{figure*}

To keep the size of the population constant, the plants with ranks $> NP$ after sorting, are eliminated. Note that in MPPA (Algorithm 3), the number of runners per plant is fixed. For this reason we refer to $n_r$ instead of $n_r^i$.

\begin{algorithm}[H]\label{MPPA}
\caption{\textbf{Pseudo-code of the Modified Plant Propagation Algorithm}}
\label{alg:MPPA}
\begin{algorithmic}[1]
\State \begin{varwidth}[t]{\linewidth} Fix the plant population size as $NP$, maximum number \par \hspace{-5mm}\hskip \algorithmicindent of function evaluations, $max\_eval$,\end{varwidth}
\State  Maximum number of generations, $max\_gen$,
\State The number of function evaluations so far, $n\_eval$,
\State \begin{varwidth}[t]{\linewidth} Create an initial population of plants \par
\hspace{-5.6mm} \hskip \algorithmicindent $pop=\{X_i\mid i=1,2,...,NP\}$,\end{varwidth}
\State Evaluate the population,
\State Set number of runners, $n_r=5$,
\State $ngen=1$, $n\_eval=NP$,
\While {$ngen$ $<$ $max\_gen$ and $n\_eval < max\_eval$}
    \State Create $\Phi$;
   \For {$i=1$ to $NP$}
      \State $\Phi_{i}=X_i$;
     \For {$k=1$ to $n_r$}
     \State \begin{varwidth}[t]{\linewidth}
     Generate a new solution $Y$ according to \par
     \hspace{-6mm} \hskip \algorithmicindent Equation \ref{eq:perturbation1};
     \end{varwidth}
     \State Evaluate it and store it in $\Phi$;
   \State Calculate $diff=F(Y)-F(X_{i})$;
     \If {diff$ \geq 0$}
     \State \begin{varwidth}[t]{\linewidth}Generate a new solution Y according \par \hspace{-4.59mm}\hskip \algorithmicindent to Equation \ref{eq:perturbation2}; \end{varwidth}
     \State Evaluate it and store it in $\Phi$;
     \State Compute $diff=f(Y)-f(X_{i})$;
    \If {diff$ \geq $0}
        \State \begin{varwidth}[t]{\linewidth} Generate a new runner using \par \hspace{-4.8mm}\hskip \algorithmicindent Equation \ref{eq:perturbation3};\end{varwidth}
        \State Evaluate it and store it in $\Phi$;
        \EndIf
      \EndIf
      \EndFor
 \EndFor
 \State Evaluate and add $\Phi$ to current population;
 \State \begin{varwidth}[t]{\linewidth} Sort the population in ascending order of the \par \hspace{-4.7mm}\hskip \algorithmicindent objective values; \end{varwidth}
 \State Update current best;
 \State Update $n\_eval$;
 \EndWhile
    \end{algorithmic}
\end{algorithm}

\section{Experimental Studies}

MPPA has been applied to a set of 18 benchmark functions of dimensions $D=30, 60 \mbox{ and } 100$, as shown in Tables 1-3. The set of experiments are carried out using $n*5000$ function evaluations, where $n$ is the dimension of the given test problem. The population size is $NP=75$, as used in \cite{gao2011modified}. We try to match the number of evaluations carried out in tests in \cite{gao2011modified}. Note that although MPPA generates more points per iteration than MABC, it runs for less generations. In this way the number of function evaluation is kept the same for both algorithms. One, therefore, can talk of a fair comparison although some readers may find the concept of ``fair comparison" hard to achieve in empirical studies on algorithms. The solution quality is listed in terms of best, worst, median, mean and standard deviation of the objective values found by each algorithm over 30 independent runs, as shown in Tables 4-6.

For functions 1, 2, 3, 5, 7, through 13, 16, and 17, MPPA found the optimum solutions while MABC did not.
For functions 4 both algorithms found the optimum.
For functions 6, 15, and 18 MPPA generated better solutions in terms of quality than MABC, although these solutions are suboptimal.
Only for function 14 did MPPA generate a solution of lower quality  than that found by MABC.

It is also important to note that MPPA outperformed ABC on all functions except Function 14, in dimension 30. Note that, at least for the time being, it is not necessary to compare with other algorithms such as the Genetic Algorithm \cite{holland1975adaption,goldberg1989genetic}, Particle Swarm Optimisation \cite{dorigo2006ant}, Differential Evolution \cite{rahnamayan2008opposition}, Harmony Search algorithm \cite{geem2001new} and others, since all of these have been outperformed by ABC as reported in \cite{karaboga2009comparative, karaboga2008performance}.

 \begin{table*}[hbt]

  \caption{Unimodal (U) and Separable (S) Test functions}
  \hspace{9mm}
    \begin{tabular}{llllll}
    \toprule
    Ftn No.    & Range & Dim         & Function & Formulation & Min \\
    \midrule
    1     & [-100, 100] & [30,60,100]     & Sphere &   $f(x)=\sum_{i=1}^{n}x^{2}_{i}$    & 0 \\[1ex]
    2     & [-100, 100] & [30,60,100]    & Elliptic  & $f(x)=\sum_{i=1}^{n}(10^{6})^{\frac{(i-1)}{(n-1)}}x^{2}_{i} $      & 0 \\[1ex]
    3     & [-10, 10] & [30,60,100]     & Axis Parallel Hyperellipsoid &   $f(x)=\sum_{i=1}^{n}|x_{i}|^{(i+1)} $     & 0 \\[1ex]
    4     & [-100, 100] & [30,60,100]     & Step  & $f(x)=\sum_{i=1}^{n}(\lfloor x_{i}+0.5\rfloor)^{2}$      & 0 \\[1ex]
    5     & [-1.28, 1.28] & [30,60,100]    & De Jong's 4 (no noise) & $f(x)=\sum_{i=1}^{n}ix^{4}_{i}$      & 0 \\[1ex]
    6     & [-1.28, 1.28] & [30,60,100]     & Quartic (noise) & $f(x)=\sum_{i=1}^{n}ix^{4}_{i}+random[0,1)$      & 0 \\[1ex]

    \bottomrule
    \end{tabular}%
  \label{tab:tab01}%
\end{table*}%
{
 \begin{table*}[htbp]

\centering
  \caption{Unimodal (U) and Non-separable (N) Test functions}
    \begin{tabular}{llllll}
    \toprule
    Ftn No.    & Range & Dim         & Function & Formulation & Min \\
    \midrule
    7     & [-10, 10] & [30,60,100]    & Sum of Different Powers &  $f(x)=\sum_{i=1}^{n}ix^{2}_{i} $     & 0 \\[1ex]
    8     & [-10, 10] & [30,60,100]    & Schwefel's Problem 2.22 & $f(x)=\sum_{i=1}^{n}|x_{i}|+\prod_{i=1}^{n}|x_{i}| $      & 0 \\[1ex]
    9     & [-100, 100] & [30,60,100]    & Schwefel's Problem 2.21 & $f(x)=max_{i}\{|x_{i}|,1\leq i\leq n\}$      & 0 \\[1ex]
    10    & [-10, 10] & [30,60,100]    & Rosenbrock &  $f(x)=\sum_{i=1}^{n-1}[100(x_{i+1}-x^{2}_{i})^{2}+(x_{i}-1)^{2}]$     & 0 \\[1ex]
    \bottomrule
    \end{tabular}%
  \label{tab:tab01}%
\end{table*}%

\fontsize{12}{12}
\selectfont
 \begin{table*}[htbp]
 \scriptsize
\centering
  \caption{Multimodal (M), Separable(S)/Non-separable(N) Test functions}
    \begin{tabular}{lllllll}
    \toprule
    Ftn No.    & Range & D     & Type     & Function & Formulation & Min \\
    \midrule
    11    & [-5.12, 5.12] & [30,60,100] & MS    & Rastrigin & $f(x)=[x^{2}_{i}-10\cos(2\pi x_{i})+10]$      & 0 \\[1ex]
    12    & [-5.12, 5.12] & [30,60,100] & M     & Non-Continuous &  $f(x)=[y^{2}_{i}-10\cos(2\pi y_{i})+10]$     & 0 \\[1ex]
        &  &  &    & Rastrigin  &  $y_{i}=\begin{cases}x_{i} & |x_{i}|< \frac{1}{2}\\ \frac{round(2x_{i})}{2}; &|x_{i}|\geq \frac{1}{2}\end{cases}$     &  \\[1ex]
    13    & [-600, 600] & [30,60,100] & MN    & Griewank & $f(x)=\frac{1}{4000}\sum_{i=1}^{n}x^{2}_{i}-\prod_{i=1}^{n}\cos(\frac{x_{i}}{\sqrt{i}})+1$      & 0 \\[1ex]
    14    & [-500, 500] & [30,60,100] & MS    & Schwefel &$f(x)=418.98288727243369*n-\sum_{i=1}^{n}x_{i}\sin(\sqrt{|x_{i}|}) $       & 0 \\[1ex]
    15    & [-32, 32] & [30,60,100] & MN    & Ackeley's Path &$f(x)=-20\exp(-0.2\sqrt{\frac{1}{n}\sum_{i=1}^{n}x^{2}_{i}})-\exp(\frac{1}{n}\sum_{i=1}^{n}\cos(2\pi x_{i}))+20+e$       & 0 \\[1ex]
    16    & [-10, 10] & [30,60,100] & M     & Alpine & $f(x)=\sum_{i=1}^{n}\mid x_{i} \cdot\sin(x_{i})+0.1 \cdot x_{i}\mid$      & 0 \\[1ex]
    17    & [-0.5, 0.5] & [30,60,100] & M     & Weierstrass & $f(x)=\sum_{i=1}^{n}(\sum_{k=0}^{k_{max}}[a^{k}\cos(2\pi b^{k}(x_{i}+0.5))]) $      &  \\[1ex]
    &  &  &    &  &  $-D\sum_{k=0}^{k_{max}}[a^{k}\cos(2\pi b^{k}(0.5)],a=0.5, b=3,k_{max}=20. $     & 0 \\[1ex]
    18    & [-100, 100] & [30,60,100] & MN    & Schaffer &$f(x)=0.5+\frac{\sin^{2}(\sqrt{\sum_{i=1}^{n}x^{2}_{i}})-0.5}{(1+0.001(\sum_{i=1}^{n}x^{2}_{i}))^{2}}$       & 0 \\[1ex]
    \bottomrule
    \end{tabular}%
  \label{tab:tab01}%
\end{table*}%

\begin{table*}[hbt]
\centering
\scriptsize
  \caption{Summary of Results Obtained with the PPA \cite{Salhi2010PPA} and MPPA Methods. Experiments were Repeated 10 Times for Each Problem and Solutions Obtained were Compared with PPA \cite{Salhi2010PPA}. a and b are the Lower and Upper Bounds On X, the Decision Variables.
}
\hspace{9mm}
    \begin{tabular}{llllccccc}
    \toprule
    Function   & Range & Dim         &  PPA & MPPA& True value& Error$_{PPA}$&Error$_{MPPA}$ \\
    \midrule
    Six hump camel back    & [-3, 2] & 2     &  -1.0316    & -1.0316 &-1.0316 &0& 0\\[1ex]
    Branin     & [-5, 15] & 2    &    0.3980    & 0.3980 &0.397887 &1.13e-04&1.13e-04\\[1ex]
    Easom    & [-100, 100] & 2     &     -0.9997     & -0.9997&-1 &3.0e-04 &3.0e-04\\[1ex]
   Goldstein Price     & [-2, 2] & 2     &    3.00536     & 3& 3&5.36e-03 &0\\[1ex]
    Martin Gaddy     & [-20, 20] & 2   &    4.4749e-08    & 0&0 &4.4749e-08 &0\\[1ex]
    Rastrigin     & [-10, 10] & 2     &  9.2e-03      & 0& 0&9.2e-03 &0\\[1ex]
    Rosenbrock     & [-5, 10] & 2    &  2.0e-04     & 0& 0&2.0e-04 &0\\[1ex]
    Schwefel     & [-500, 500] & 2     &   -833.203     & -837.9603& -837.9658 &4.76 &5.5e-03\\[1ex]
    \bottomrule
    \end{tabular}%
  \label{tab:tab01}%
\end{table*}%

\begin{table*}[htbp]
 \scriptsize
 \centering
  \caption{Results obtained by MPPA, ABC and MABC on test functions of Table 1.}
    \begin{tabular}{clrrrrrr}
    \toprule
    Fun   & Dim   & Algorithm       & Best  & Worst & Median & Mean  & SD \\
    \midrule
    1     & 30    & ABC   & \multicolumn{1}{l}{2.02e - 10} & \multicolumn{1}{l}{1.07e - 09} & \multicolumn{1}{l}{2.02e - 10} & \multicolumn{1}{l}{5.21e - 10} & \multicolumn{1}{l}{2.46e - 10} \\
          &       & MABC  & \multicolumn{1}{l}{1.69e - 32} & \multicolumn{1}{l}{2.48e - 31} & \multicolumn{1}{l}{1.80e - 32} & \multicolumn{1}{l}{9.43e - 32} & \multicolumn{1}{l}{6.67e - 32} \\
          &       & MPPA  & \multicolumn{1}{l}{0} & \multicolumn{1}{l}{0} & \multicolumn{1}{l}{0} & \multicolumn{1}{l}{0} & \multicolumn{1}{l}{0} \\
          & 60    & ABC   & \multicolumn{1}{l}{3.49e - 10} & \multicolumn{1}{l}{3.27e - 09} & \multicolumn{1}{l}{3.27e - 09} & \multicolumn{1}{l}{1.09e - 09} & \multicolumn{1}{l}{9.37e - 10} \\
          &       & MABC  & \multicolumn{1}{l}{9.29e - 30} & \multicolumn{1}{l}{1.55e - 28} & \multicolumn{1}{l}{4.68e - 29} & \multicolumn{1}{l}{6.03e - 29} & \multicolumn{1}{l}{4.31e - 29} \\
          &       & MPPA  & \multicolumn{1}{l}{0} & \multicolumn{1}{l}{0} & \multicolumn{1}{l}{0} & \multicolumn{1}{l}{0} & \multicolumn{1}{l}{0} \\
          & 100   & ABC   & \multicolumn{1}{l}{3.87e - 10} & \multicolumn{1}{l}{3.11e - 09} & \multicolumn{1}{l}{3.87e - 10} & \multicolumn{1}{l}{1.64e - 09} & \multicolumn{1}{l}{9.85e - 10} \\
          &       & MABC  & \multicolumn{1}{l}{4.65e - 28} & \multicolumn{1}{l}{2.63e - 27} & \multicolumn{1}{l}{1.86e - 27} & \multicolumn{1}{l}{1.43e - 27} & \multicolumn{1}{l}{8.12e - 28} \\
          &       & MPPA  & \multicolumn{1}{l}{0} & \multicolumn{1}{l}{0} & \multicolumn{1}{l}{0} & \multicolumn{1}{l}{0} & \multicolumn{1}{l}{0} \\
    2     & 30    & ABC   & \multicolumn{1}{l}{2.65e - 07} & \multicolumn{1}{l}{1.32e - 05} & \multicolumn{1}{l}{7.71e - 07} & \multicolumn{1}{l}{4.10e - 06} & \multicolumn{1}{l}{3.85e - 06} \\
          &       & MABC  & \multicolumn{1}{l}{2.55e - 29} & \multicolumn{1}{l}{1.99e - 27} & \multicolumn{1}{l}{1.99e - 27} & \multicolumn{1}{l}{3.66e - 28} & \multicolumn{1}{l}{5.96e - 28} \\
          &       & MPPA  & \multicolumn{1}{l}{0} & \multicolumn{1}{l}{0} & \multicolumn{1}{l}{0} & \multicolumn{1}{l}{0} & \multicolumn{1}{l}{0} \\
          & 60    & ABC   & \multicolumn{1}{l}{2.14e - 07} & \multicolumn{1}{l}{6.25e - 06} & \multicolumn{1}{l}{4.53e - 06} & \multicolumn{1}{l}{2.31e - 06} & \multicolumn{1}{l}{2.18e - 06} \\
          &       & MABC  & \multicolumn{1}{l}{3.65e - 26} & \multicolumn{1}{l}{8.69e - 25} & \multicolumn{1}{l}{3.65e - 26} & \multicolumn{1}{l}{3.51e - 25} & \multicolumn{1}{l}{2.72e - 25} \\
          &       & MPPA  & \multicolumn{1}{l}{0} & \multicolumn{1}{l}{0} & \multicolumn{1}{l}{0} & \multicolumn{1}{l}{0} & \multicolumn{1}{l}{0} \\
          & 100   & ABC   & \multicolumn{1}{l}{4.00e - 07} & \multicolumn{1}{l}{6.46e - 06} & \multicolumn{1}{l}{1.37e - 06} & \multicolumn{1}{l}{1.79e - 06} & \multicolumn{1}{l}{1.63e - 06} \\
          &       & MABC  & \multicolumn{1}{l}{3.17e - 24} & \multicolumn{1}{l}{3.73e - 24} & \multicolumn{1}{l}{3.17e - 24} & \multicolumn{1}{l}{3.52e - 24} & \multicolumn{1}{l}{2.47e - 25} \\
          &       & MPPA  & \multicolumn{1}{l}{0} & \multicolumn{1}{l}{0} & \multicolumn{1}{l}{0} & \multicolumn{1}{l}{0} & \multicolumn{1}{l}{0} \\
    3     & 30    & ABC   & \multicolumn{1}{l}{9.03e - 12} & \multicolumn{1}{l}{4.62e - 11} & \multicolumn{1}{l}{1.56e - 11} & \multicolumn{1}{l}{2.22e - 11} & \multicolumn{1}{l}{1.14e - 11} \\
          &       & MABC  & \multicolumn{1}{l}{3.57e - 33} & \multicolumn{1}{l}{5.29e - 32} & \multicolumn{1}{l}{4.62e - 33} & \multicolumn{1}{l}{2.10e - 32} & \multicolumn{1}{l}{1.56e - 32} \\
          &       & MPPA  & \multicolumn{1}{l}{0} & \multicolumn{1}{l}{0} & \multicolumn{1}{l}{0} & \multicolumn{1}{l}{0} & \multicolumn{1}{l}{0} \\
          & 60    & ABC   & \multicolumn{1}{l}{7.95e - 11} & \multicolumn{1}{l}{3.17e - 10} & \multicolumn{1}{l}{3.17e - 10} & \multicolumn{1}{l}{1.89e - 10} & \multicolumn{1}{l}{9.14e - 11} \\
          &       & MABC  & \multicolumn{1}{l}{5.61e - 30} & \multicolumn{1}{l}{3.29e - 29} & \multicolumn{1}{l}{6.24e - 30} & \multicolumn{1}{l}{1.39e - 29} & \multicolumn{1}{l}{8.84e - 30} \\
          &       & MPPA  & \multicolumn{1}{l}{0} & \multicolumn{1}{l}{0} & \multicolumn{1}{l}{0} & \multicolumn{1}{l}{0} & \multicolumn{1}{l}{0} \\
          & 100   & ABC   & \multicolumn{1}{l}{2.82e - 10} & \multicolumn{1}{l}{2.85e - 09} & \multicolumn{1}{l}{4.21e - 10} & \multicolumn{1}{l}{1.25e - 09} & \multicolumn{1}{l}{9.75e - 10} \\
          &       & MABC  & \multicolumn{1}{l}{1.50e - 28} & \multicolumn{1}{l}{8.74e - 28} & \multicolumn{1}{l}{1.72e - 28} & \multicolumn{1}{l}{4.46e - 28} & \multicolumn{1}{l}{2.08e - 28} \\
          &       & MPPA  & \multicolumn{1}{l}{0} & \multicolumn{1}{l}{0} & \multicolumn{1}{l}{0} & \multicolumn{1}{l}{0} & \multicolumn{1}{l}{0} \\

    4     & 30    & ABC   & \multicolumn{1}{l}{0} & \multicolumn{1}{l}{0} & \multicolumn{1}{l}{0} & \multicolumn{1}{l}{0} & \multicolumn{1}{l}{0} \\
          &       & MABC  & \multicolumn{1}{l}{0} & \multicolumn{1}{l}{0} & \multicolumn{1}{l}{0} & \multicolumn{1}{l}{0} & \multicolumn{1}{l}{0} \\
          &       & MPPA  & \multicolumn{1}{l}{0} & \multicolumn{1}{l}{0} & \multicolumn{1}{l}{0} & \multicolumn{1}{l}{0} & \multicolumn{1}{l}{0} \\
          & 60    & ABC   & \multicolumn{1}{l}{0} & \multicolumn{1}{l}{0} & \multicolumn{1}{l}{0} & \multicolumn{1}{l}{0} & \multicolumn{1}{l}{0} \\
          &       & MABC  & \multicolumn{1}{l}{0} & \multicolumn{1}{l}{0} & \multicolumn{1}{l}{0} & \multicolumn{1}{l}{0} & \multicolumn{1}{l}{0} \\
          &       & MPPA  & \multicolumn{1}{l}{0} & \multicolumn{1}{l}{0} & \multicolumn{1}{l}{0} & \multicolumn{1}{l}{0} & \multicolumn{1}{l}{0} \\
          & 100   & ABC   & \multicolumn{1}{l}{0} & \multicolumn{1}{l}{0} & \multicolumn{1}{l}{0} & \multicolumn{1}{l}{0} & \multicolumn{1}{l}{0} \\
          &       & MABC  & \multicolumn{1}{l}{0} & \multicolumn{1}{l}{0} & \multicolumn{1}{l}{0} & \multicolumn{1}{l}{0} & \multicolumn{1}{l}{0} \\
          &       & MPPA  & \multicolumn{1}{l}{0} & \multicolumn{1}{l}{0} & \multicolumn{1}{l}{0} & \multicolumn{1}{l}{0} & \multicolumn{1}{l}{0} \\
                5& 30    & ABC   & \multicolumn{1}{l}{1.26e - 29} & \multicolumn{1}{l}{1.86e - 28} & \multicolumn{1}{l}{1.64e - 29} & \multicolumn{1}{l}{5.51e - 29} & \multicolumn{1}{l}{6.70e - 29} \\
          &       & MABC  & \multicolumn{1}{l}{8.80e - 69} & \multicolumn{1}{l}{5.97e - 67} & \multicolumn{1}{l}{1.01e - 68} & \multicolumn{1}{l}{1.45e - 67} & \multicolumn{1}{l}{2.28e - 67} \\
          &       & MPPA  & \multicolumn{1}{l}{0} & \multicolumn{1}{l}{0} & \multicolumn{1}{l}{0} & \multicolumn{1}{l}{0} & \multicolumn{1}{l}{0} \\
          & 60    & ABC   & \multicolumn{1}{l}{4.90e - 28} & \multicolumn{1}{l}{2.00e - 26} & \multicolumn{1}{l}{4.90e - 28} & \multicolumn{1}{l}{6.53e - 27} & \multicolumn{1}{l}{7.23e - 27} \\
          &       & MABC  & \multicolumn{1}{l}{4.49e - 64} & \multicolumn{1}{l}{2.37e - 61} & \multicolumn{1}{l}{2.37e - 61} & \multicolumn{1}{l}{5.00e - 62} & \multicolumn{1}{l}{9.38e - 62} \\
          &       & MPPA  & \multicolumn{1}{l}{0} & \multicolumn{1}{l}{0} & \multicolumn{1}{l}{0} & \multicolumn{1}{l}{0} & \multicolumn{1}{l}{0} \\
          & 100   & ABC   & \multicolumn{1}{l}{1.31e - 26} & \multicolumn{1}{l}{1.41e - 25} & \multicolumn{1}{l}{1.34e - 26} & \multicolumn{1}{l}{5.65e - 26} & \multicolumn{1}{l}{4.90e - 26} \\
          &       & MABC  & \multicolumn{1}{l}{3.99e - 61} & \multicolumn{1}{l}{1.56e - 59} & \multicolumn{1}{l}{3.99e - 61} & \multicolumn{1}{l}{5.72e - 60} & \multicolumn{1}{l}{5.32e - 60} \\
          &       & MPPA  & \multicolumn{1}{l}{0} & \multicolumn{1}{l}{0} & \multicolumn{1}{l}{0} & \multicolumn{1}{l}{0} & \multicolumn{1}{l}{0} \\
         6& 30    & ABC   & \multicolumn{1}{l}{6.03e - 02} & \multicolumn{1}{l}{1.27e - 01} & \multicolumn{1}{l}{7.67e - 02} & \multicolumn{1}{l}{8.74e - 02} & \multicolumn{1}{l}{1.77e - 02} \\
          &       & MABC  & \multicolumn{1}{l}{1.84e - 02} & \multicolumn{1}{l}{4.58e - 02} & \multicolumn{1}{l}{4.48e - 02} & \multicolumn{1}{l}{3.71e - 02} & \multicolumn{1}{l}{8.53e - 03} \\
          &       & MPPA  & \multicolumn{1}{l}{1.98e - 06} & \multicolumn{1}{l}{5.46e - 05} & \multicolumn{1}{l}{1.41e - 05} & \multicolumn{1}{l}{1.78e - 05} & \multicolumn{1}{l}{1.39e - 05} \\
          & 60    & ABC   & \multicolumn{1}{l}{1.87e - 01} & \multicolumn{1}{l}{2.65e - 01} & \multicolumn{1}{l}{2.43e - 01} & \multicolumn{1}{l}{2.39e - 01} & \multicolumn{1}{l}{2.86e - 02} \\
          &       & MABC  & \multicolumn{1}{l}{9.20e - 02} & \multicolumn{1}{l}{1.33e - 01} & \multicolumn{1}{l}{1.21e - 01} & \multicolumn{1}{l}{1.14e - 01} & \multicolumn{1}{l}{1.16e - 02} \\
          &       & MPPA  & \multicolumn{1}{l}{2.78e - 08} & \multicolumn{1}{l}{1.54e - 05} & \multicolumn{1}{l}{4.25e - 06} & \multicolumn{1}{l}{5.12e - 06} & \multicolumn{1}{l}{4.25e - 06} \\
          & 100   & ABC   & \multicolumn{1}{l}{3.96e - 01} & \multicolumn{1}{l}{4.92e - 01} & \multicolumn{1}{l}{4.70e - 01} & \multicolumn{1}{l}{4.55e - 01} & \multicolumn{1}{l}{3.20e - 02} \\
          &       & MABC  & \multicolumn{1}{l}{1.64e - 01} & \multicolumn{1}{l}{2.56e - 01} & \multicolumn{1}{l}{2.56e - 01} & \multicolumn{1}{l}{2.31e - 01} & \multicolumn{1}{l}{2.79e - 02} \\
          &       & MPPA  & \multicolumn{1}{l}{7.97e - 08} & \multicolumn{1}{l}{2.44e - 05} & \multicolumn{1}{l}{3.28e - 06} & \multicolumn{1}{l}{5.28e - 06} & \multicolumn{1}{l}{5.80e - 06} \\
       \bottomrule
    \end{tabular}%
  \label{tab:addlabel}%
\end{table*}%

\begin{table*}[htbp]
 \scriptsize
  \centering
  \caption{Results obtained by MPPA, ABC and MABC on test functions of Table 2}
    \begin{tabular}{clrrrrrr}
    \toprule
    Fun   & Dim   & Algorithm      & Best  & Worst & Median & Mean  & SD \\
    \midrule
          7& 30    & ABC   & \multicolumn{1}{l}{1.34e - 18} & \multicolumn{1}{l}{4.29e - 16} & \multicolumn{1}{l}{1.82e - 16} & \multicolumn{1}{l}{1.45e - 16} & \multicolumn{1}{l}{1.55e - 16} \\
          &       & MABC  & \multicolumn{1}{l}{1.26e - 74} & \multicolumn{1}{l}{1.34e - 68} & \multicolumn{1}{l}{7.02e - 71} & \multicolumn{1}{l}{2.70e - 69} & \multicolumn{1}{l}{5.38e - 69} \\
          &       & MPPA  & \multicolumn{1}{l}{0} & \multicolumn{1}{l}{0} & \multicolumn{1}{l}{0} & \multicolumn{1}{l}{0} & \multicolumn{1}{l}{0} \\
          & 60    & ABC   & \multicolumn{1}{l}{4.88e - 11} & \multicolumn{1}{l}{7.62e - 10} & \multicolumn{1}{l}{4.88e - 11} & \multicolumn{1}{l}{2.14e - 10} & \multicolumn{1}{l}{2.75e - 10} \\
          &       & MABC  & \multicolumn{1}{l}{5.16e - 65} & \multicolumn{1}{l}{9.54e - 62} & \multicolumn{1}{l}{3.59e - 64} & \multicolumn{1}{l}{3.00e - 62} & \multicolumn{1}{l}{3.87e - 62} \\
          &       & MPPA  & \multicolumn{1}{l}{0} & \multicolumn{1}{l}{0} & \multicolumn{1}{l}{0} & \multicolumn{1}{l}{0} & \multicolumn{1}{l}{0} \\
          & 100   & ABC   & \multicolumn{1}{l}{2.52e - 08} & \multicolumn{1}{l}{2.05e - 06} & \multicolumn{1}{l}{6.75e - 08} & \multicolumn{1}{l}{4.83e - 07} & \multicolumn{1}{l}{7.88e - 07} \\
          &       & MABC  & \multicolumn{1}{l}{7.76e - 52} & \multicolumn{1}{l}{8.74e - 48} & \multicolumn{1}{l}{8.01e - 52} & \multicolumn{1}{l}{1.92e - 48} & \multicolumn{1}{l}{3.42e - 48} \\
          &       & MPPA  & \multicolumn{1}{l}{0} & \multicolumn{1}{l}{0} & \multicolumn{1}{l}{0} & \multicolumn{1}{l}{0} & \multicolumn{1}{l}{0} \\
         8& 30    & ABC   & \multicolumn{1}{l}{1.28e - 06} & \multicolumn{1}{l}{2.63e - 06} & \multicolumn{1}{l}{1.28e - 06} & \multicolumn{1}{l}{1.83e - 06} & \multicolumn{1}{l}{4.80e - 07} \\
          &       & MABC  & \multicolumn{1}{l}{8.41e - 18} & \multicolumn{1}{l}{4.09e - 17} & \multicolumn{1}{l}{8.41e - 18} & \multicolumn{1}{l}{2.40e - 17} & \multicolumn{1}{l}{9.02e - 18} \\
          &       & MPPA  & \multicolumn{1}{l}{0} & \multicolumn{1}{l}{0} & \multicolumn{1}{l}{0} & \multicolumn{1}{l}{0} & \multicolumn{1}{l}{0} \\
          & 60    & ABC   & \multicolumn{1}{l}{5.77e - 06} & \multicolumn{1}{l}{9.29e - 06} & \multicolumn{1}{l}{6.81e - 06} & \multicolumn{1}{l}{7.23e - 06} & \multicolumn{1}{l}{1.28e - 06} \\
          &       & MABC  & \multicolumn{1}{l}{5.30e - 16} & \multicolumn{1}{l}{8.72e - 16} & \multicolumn{1}{l}{8.33e - 16} & \multicolumn{1}{l}{6.96e - 16} & \multicolumn{1}{l}{1.20e - 16} \\
          &       & MPPA  & \multicolumn{1}{l}{0} & \multicolumn{1}{l}{0} & \multicolumn{1}{l}{0} & \multicolumn{1}{l}{0} & \multicolumn{1}{l}{0} \\
          & 100   & ABC   & \multicolumn{1}{l}{1.28e - 06} & \multicolumn{1}{l}{1.59e - 05} & \multicolumn{1}{l}{1.29e - 05} & \multicolumn{1}{l}{1.30e - 05} & \multicolumn{1}{l}{1.93e - 06} \\
          &       & MABC  & \multicolumn{1}{l}{2.08e - 15} & \multicolumn{1}{l}{6.63e - 15} & \multicolumn{1}{l}{6.46e - 15} & \multicolumn{1}{l}{4.41e - 15} & \multicolumn{1}{l}{1.50e - 15} \\
          &       & MPPA  & \multicolumn{1}{l}{0} & \multicolumn{1}{l}{0} & \multicolumn{1}{l}{0} & \multicolumn{1}{l}{0} & \multicolumn{1}{l}{0} \\
         9& 30    & ABC   & \multicolumn{1}{l}{1.30e + 01} & \multicolumn{1}{l}{2.07e + 01} & \multicolumn{1}{l}{1.52e + 01} & \multicolumn{1}{l}{1.80e + 01} & \multicolumn{1}{l}{2.25e - 00} \\
          &       & MABC  & \multicolumn{1}{l}{9.16e - 00} & \multicolumn{1}{l}{1.42e + 01} & \multicolumn{1}{l}{1.26e + 01} & \multicolumn{1}{l}{1.02e + 01} & \multicolumn{1}{l}{1.49e - 00} \\
          &       & MPPA  & \multicolumn{1}{l}{0} & \multicolumn{1}{l}{0} & \multicolumn{1}{l}{0} & \multicolumn{1}{l}{0} & \multicolumn{1}{l}{0} \\
          & 60    & ABC   & \multicolumn{1}{l}{3.79e + 01} & \multicolumn{1}{l}{4.65e + 01} & \multicolumn{1}{l}{4.18e + 01} & \multicolumn{1}{l}{4.22e + 01} & \multicolumn{1}{l}{2.73e - 00} \\
          &       & MABC  & \multicolumn{1}{l}{2.92e + 01} & \multicolumn{1}{l}{4.03e + 01} & \multicolumn{1}{l}{3.80e + 01} & \multicolumn{1}{l}{3.77e + 01} & \multicolumn{1}{l}{3.14e - 00} \\
          &       & MPPA  & \multicolumn{1}{l}{0} & \multicolumn{1}{l}{0} & \multicolumn{1}{l}{0} & \multicolumn{1}{l}{0} & \multicolumn{1}{l}{0} \\
          & 100   & ABC   & \multicolumn{1}{l}{5.30e + 01} & \multicolumn{1}{l}{6.14e + 01} & \multicolumn{1}{l}{5.30e + 01} & \multicolumn{1}{l}{5.76e + 01} & \multicolumn{1}{l}{2.74e - 00} \\
          &       & MABC  & \multicolumn{1}{l}{5.70e + 01} & \multicolumn{1}{l}{6.23e + 01} & \multicolumn{1}{l}{5.76e + 01} & \multicolumn{1}{l}{5.98e + 01} & \multicolumn{1}{l}{1.60e - 00} \\
    10    & 30    & ABC   & \multicolumn{1}{l}{2.12e - 02} & \multicolumn{1}{l}{2.20e - 00} & \multicolumn{1}{l}{5.56e - 01} & \multicolumn{1}{l}{4.23e - 01} & \multicolumn{1}{l}{4.34e - 01} \\
          &       & MABC  & \multicolumn{1}{l}{4.09e - 02} & \multicolumn{1}{l}{1.95e - 00} & \multicolumn{1}{l}{2.34e - 01} & \multicolumn{1}{l}{6.11e - 01} & \multicolumn{1}{l}{4.55e - 01} \\
          &       & MPPA  & \multicolumn{1}{l}{0} & \multicolumn{1}{l}{8.24e - 03} & \multicolumn{1}{l}{1.42e - 05} & \multicolumn{1}{l}{4.24e - 04} & \multicolumn{1}{l}{1.60e - 03} \\
          & 60    & ABC   & \multicolumn{1}{l}{1.88e - 01} & \multicolumn{1}{l}{5.75e - 00} & \multicolumn{1}{l}{1.80e - 00} & \multicolumn{1}{l}{1.86e - 00} & \multicolumn{1}{l}{1.36e - 00} \\
          &       & MABC  & \multicolumn{1}{l}{2.17e - 01} & \multicolumn{1}{l}{5.26e - 00} & \multicolumn{1}{l}{1.30e - 00} & \multicolumn{1}{l}{1.51e - 00} & \multicolumn{1}{l}{1.34e - 00} \\
          &       & MPPA  & \multicolumn{1}{l}{0} & \multicolumn{1}{l}{2.13e - 03} & \multicolumn{1}{l}{4.02e - 05} & \multicolumn{1}{l}{1.83e - 04} & \multicolumn{1}{l}{5.08e - 04} \\
          & 100   & ABC   & \multicolumn{1}{l}{5.48e - 01} & \multicolumn{1}{l}{3.94e - 00} & \multicolumn{1}{l}{6.33e - 01} & \multicolumn{1}{l}{1.59e - 00} & \multicolumn{1}{l}{1.23e - 00} \\
          &       & MABC  & \multicolumn{1}{l}{5.13e - 01} & \multicolumn{1}{l}{7.77e - 00} & \multicolumn{1}{l}{6.71e - 01} & \multicolumn{1}{l}{1.98e - 00} & \multicolumn{1}{l}{1.30e - 00} \\
          &       & MPPA  & \multicolumn{1}{l}{0} & \multicolumn{1}{l}{1.22e - 03} & \multicolumn{1}{l}{1.82e - 04} & \multicolumn{1}{l}{3.52e - 04} & \multicolumn{1}{l}{3.81e - 04} \\
       \bottomrule
    \end{tabular}%
  \label{tab:addlabel}%
\end{table*}%

\begin{table*}[htbp]
 \scriptsize
  \centering
  \caption{Results obtained by MPPA, ABC and MABC on test functions of Table 3}
    \begin{tabular}{clrrrrrr}
    \toprule
    Fun   & Dim   & Algorithm       & Best  & Worst & Median & Mean  & SD \\
    \midrule

    11    & 30    & ABC   & \multicolumn{1}{l}{3.58e - 10} & \multicolumn{1}{l}{1.46e - 01} & \multicolumn{1}{l}{5.50e - 10} & \multicolumn{1}{l}{4.81e - 03} & \multicolumn{1}{l}{2.57e - 02} \\
          &       & MABC  & \multicolumn{1}{l}{0} & \multicolumn{1}{l}{0} & \multicolumn{1}{l}{0} & \multicolumn{1}{l}{0} & \multicolumn{1}{l}{0} \\
          &       & MPPA  & \multicolumn{1}{l}{0} & \multicolumn{1}{l}{0} & \multicolumn{1}{l}{0} & \multicolumn{1}{l}{0} & \multicolumn{1}{l}{0} \\
          & 60    & ABC   & \multicolumn{1}{l}{3.21e - 10} & \multicolumn{1}{l}{1.99e - 00} & \multicolumn{1}{l}{1.28e - 05} & \multicolumn{1}{l}{3.71e - 01} & \multicolumn{1}{l}{5.97e - 01} \\
          &       & MABC  & \multicolumn{1}{l}{0} & \multicolumn{1}{l}{0} & \multicolumn{1}{l}{0} & \multicolumn{1}{l}{0} & \multicolumn{1}{l}{0} \\
          &       & MPPA  & \multicolumn{1}{l}{0} & \multicolumn{1}{l}{0} & \multicolumn{1}{l}{0} & \multicolumn{1}{l}{0} & \multicolumn{1}{l}{0} \\
          & 100   & ABC   & \multicolumn{1}{l}{5.40e - 09} & \multicolumn{1}{l}{1.99e - 00} & \multicolumn{1}{l}{1.99e - 00} & \multicolumn{1}{l}{1.10e - 00} & \multicolumn{1}{l}{8.21e - 01} \\
          &       & MABC  & \multicolumn{1}{l}{0} & \multicolumn{1}{l}{0} & \multicolumn{1}{l}{0} & \multicolumn{1}{l}{0} & \multicolumn{1}{l}{0} \\
          &       & MPPA  & \multicolumn{1}{l}{0} & \multicolumn{1}{l}{0} & \multicolumn{1}{l}{0} & \multicolumn{1}{l}{0} & \multicolumn{1}{l}{0} \\
    12    & 30    & ABC   & \multicolumn{1}{l}{8.96e - 10} & \multicolumn{1}{l}{1.00e - 00} & \multicolumn{1}{l}{1.77e - 08} & \multicolumn{1}{l}{1.12e - 01} & \multicolumn{1}{l}{2.97e - 01} \\
          &       & MABC  & \multicolumn{1}{l}{0} & \multicolumn{1}{l}{0} & \multicolumn{1}{l}{0} & \multicolumn{1}{l}{0} & \multicolumn{1}{l}{0} \\
          &       & MPPA  & \multicolumn{1}{l}{0} & \multicolumn{1}{l}{0} & \multicolumn{1}{l}{0} & \multicolumn{1}{l}{0} & \multicolumn{1}{l}{0} \\
          & 60    & ABC   & \multicolumn{1}{l}{2.17e - 07} & \multicolumn{1}{l}{3.025e - 00} & \multicolumn{1}{l}{1.09e - 00} & \multicolumn{1}{l}{1.47e - 00} & \multicolumn{1}{l}{9.47e - 01} \\
          &       & MABC  & \multicolumn{1}{l}{0} & \multicolumn{1}{l}{0} & \multicolumn{1}{l}{0} & \multicolumn{1}{l}{0} & \multicolumn{1}{l}{0} \\
          &       & MPPA  & \multicolumn{1}{l}{0} & \multicolumn{1}{l}{0} & \multicolumn{1}{l}{0} & \multicolumn{1}{l}{0} & \multicolumn{1}{l}{0} \\
          & 100   & ABC   & \multicolumn{1}{l}{2.00e - 00} & \multicolumn{1}{l}{7.21e - 00} & \multicolumn{1}{l}{5.02e - 00} & \multicolumn{1}{l}{4.74e - 00} & \multicolumn{1}{l}{2.01e - 00} \\
          &       & MABC  & \multicolumn{1}{l}{0} & \multicolumn{1}{l}{0} & \multicolumn{1}{l}{0} & \multicolumn{1}{l}{0} & \multicolumn{1}{l}{0} \\
          &       & MPPA  & \multicolumn{1}{l}{0} & \multicolumn{1}{l}{0} & \multicolumn{1}{l}{0} & \multicolumn{1}{l}{0} & \multicolumn{1}{l}{0} \\
    13    & 30    & ABC   & \multicolumn{1}{l}{4.74e - 12} & \multicolumn{1}{l}{1.35e - 07} & \multicolumn{1}{l}{1.77e - 08} & \multicolumn{1}{l}{1.61e - 08} & \multicolumn{1}{l}{3.99e - 08} \\
          &       & MABC  & \multicolumn{1}{l}{0} & \multicolumn{1}{l}{0} & \multicolumn{1}{l}{0} & \multicolumn{1}{l}{0} & \multicolumn{1}{l}{0} \\
          &       & MPPA  & \multicolumn{1}{l}{0} & \multicolumn{1}{l}{0} & \multicolumn{1}{l}{0} & \multicolumn{1}{l}{0} & \multicolumn{1}{l}{0} \\
          & 60    & ABC   & \multicolumn{1}{l}{7.99e - 12} & \multicolumn{1}{l}{1.05e - 09} & \multicolumn{1}{l}{1.45e - 11} & \multicolumn{1}{l}{1.39e - 10} & \multicolumn{1}{l}{3.10e - 10} \\
          &       & MABC  & \multicolumn{1}{l}{0} & \multicolumn{1}{l}{0} & \multicolumn{1}{l}{0} & \multicolumn{1}{l}{0} & \multicolumn{1}{l}{0} \\
          &       & MPPA  & \multicolumn{1}{l}{0} & \multicolumn{1}{l}{0} & \multicolumn{1}{l}{0} & \multicolumn{1}{l}{0} & \multicolumn{1}{l}{0} \\
          & 100   & ABC   & \multicolumn{1}{l}{6.29e - 13} & \multicolumn{1}{l}{3.05e - 08} & \multicolumn{1}{l}{1.01e - 10} & \multicolumn{1}{l}{2.01e - 09} & \multicolumn{1}{l}{1.32e - 09} \\
          &       & MABC  & \multicolumn{1}{l}{0} & \multicolumn{1}{l}{0} & \multicolumn{1}{l}{0} & \multicolumn{1}{l}{0} & \multicolumn{1}{l}{0} \\
          &       & MPPA  & \multicolumn{1}{l}{0} & \multicolumn{1}{l}{0} & \multicolumn{1}{l}{0} & \multicolumn{1}{l}{0} & \multicolumn{1}{l}{0} \\
    14    & 30    & ABC   & \multicolumn{1}{l}{1.54e - 06} & \multicolumn{1}{l}{2.37e + 02} & \multicolumn{1}{l}{3.76e - 01} & \multicolumn{1}{l}{8.86e + 01} & \multicolumn{1}{l}{8.62e + 01} \\
          &       & MABC  & \multicolumn{1}{l}{- 1.81e - 12} & \multicolumn{1}{l}{0} & \multicolumn{1}{l}{0} & \multicolumn{1}{l}{- 1.21e - 13} & \multicolumn{1}{l}{4.53e - 13} \\
          &       & MPPA  & \multicolumn{1}{l}{3.82e - 04} & \multicolumn{1}{l}{3.55e + 03} & \multicolumn{1}{l}{5.36e - 04} & \multicolumn{1}{l}{3.55e + 02} & \multicolumn{1}{l}{1.08e + 03} \\
          & 60    & ABC   & \multicolumn{1}{l}{3.55e + 02} & \multicolumn{1}{l}{7.69e + 02} & \multicolumn{1}{l}{7.69e + 02} & \multicolumn{1}{l}{5.40e + 02} & \multicolumn{1}{l}{1.41e + 02} \\
          &       & MABC  & \multicolumn{1}{l}{2.91e - 11} & \multicolumn{1}{l}{3.63e - 11} & \multicolumn{1}{l}{2.91e - 11} & \multicolumn{1}{l}{3.56e - 11} & \multicolumn{1}{l}{2.18e - 12} \\
          &       & MPPA  & \multicolumn{1}{l}{7.63e - 04} & \multicolumn{1}{l}{7.10e +03} & \multicolumn{1}{l}{7.79e - 04} & \multicolumn{1}{l}{7.10e - 02} & \multicolumn{1}{l}{2.16e +03} \\
          & 100   & ABC   & \multicolumn{1}{l}{7.81e + 02} & \multicolumn{1}{l}{1.55e + 03} & \multicolumn{1}{l}{1.51e + 03} & \multicolumn{1}{l}{1.29e + 03} & \multicolumn{1}{l}{2.23e + 02} \\
          &       & MABC  & \multicolumn{1}{l}{1.09e - 10} & \multicolumn{1}{l}{1.23e - 10} & \multicolumn{1}{l}{1.16e - 10} & \multicolumn{1}{l}{1.19e - 10} & \multicolumn{1}{l}{4.06e - 12} \\
          &       & MPPA  & \multicolumn{1}{l}{1.27e - 03} & \multicolumn{1}{l}{1.18e +04} & \multicolumn{1}{l}{1.30e -03} & \multicolumn{1}{l}{7.89e +02} & \multicolumn{1}{l}{3.00e +03} \\
          15    & 30    & ABC   & \multicolumn{1}{l}{2.26e - 06} & \multicolumn{1}{l}{8.32e - 06} & \multicolumn{1}{l}{7.17e - 06} & \multicolumn{1}{l}{4.83e - 06} & \multicolumn{1}{l}{2.12e - 06} \\
          &       & MABC  & \multicolumn{1}{l}{3.64e - 14} & \multicolumn{1}{l}{4.35e - 14} & \multicolumn{1}{l}{3.99e - 14} & \multicolumn{1}{l}{4.13e - 14} & \multicolumn{1}{l}{2.17e - 15} \\
          &       & MPPA  & \multicolumn{1}{l}{8.88e - 16} & \multicolumn{1}{l}{8.88e - 16} & \multicolumn{1}{l}{8.88e - 16} & \multicolumn{1}{l}{8.88e - 16} & \multicolumn{1}{l}{0} \\
          & 60    & ABC   & \multicolumn{1}{l}{2.44e - 06} & \multicolumn{1}{l}{1.57e - 05} & \multicolumn{1}{l}{2.44e - 06} & \multicolumn{1}{l}{7.79e - 06} & \multicolumn{1}{l}{3.63e - 06} \\
          &       & MABC  & \multicolumn{1}{l}{1.14e - 13} & \multicolumn{1}{l}{1.57e - 13} & \multicolumn{1}{l}{1.32e - 13} & \multicolumn{1}{l}{1.37e - 13} & \multicolumn{1}{l}{1.24e - 14} \\
          &       & MPPA  & \multicolumn{1}{l}{8.88e - 16} & \multicolumn{1}{l}{8.88e - 16} & \multicolumn{1}{l}{8.88e - 16} & \multicolumn{1}{l}{8.88e - 16} & \multicolumn{1}{l}{0} \\
          & 100   & ABC   & \multicolumn{1}{l}{5.12e - 06} & \multicolumn{1}{l}{1.38e - 05} & \multicolumn{1}{l}{5.94e - 06} & \multicolumn{1}{l}{1.02e - 05} & \multicolumn{1}{l}{2.92e - 06} \\
          &       & MABC  & \multicolumn{1}{l}{3.27e - 13} & \multicolumn{1}{l}{3.98e - 13} & \multicolumn{1}{l}{3.27e - 13} & \multicolumn{1}{l}{3.56e - 13} & \multicolumn{1}{l}{2.29e - 14} \\
          &       & MPPA  & \multicolumn{1}{l}{8.88e - 16} & \multicolumn{1}{l}{8.88e - 16} & \multicolumn{1}{l}{8.88e - 16} & \multicolumn{1}{l}{8.88e - 16} & \multicolumn{1}{l}{0} \\
    16    & 30    & ABC   & \multicolumn{1}{l}{2.99e - 05} & \multicolumn{1}{l}{1.05e - 04} & \multicolumn{1}{l}{1.05e - 04} & \multicolumn{1}{l}{7.66e - 05} & \multicolumn{1}{l}{2.76e - 05} \\
          &       & MABC  & \multicolumn{1}{l}{3.74e - 18} & \multicolumn{1}{l}{6.54e - 16} & \multicolumn{1}{l}{1.48e - 17} & \multicolumn{1}{l}{1.58e - 16} & \multicolumn{1}{l}{2.48e - 16} \\
          &       & MPPA  & \multicolumn{1}{l}{0} & \multicolumn{1}{l}{0} & \multicolumn{1}{l}{0} & \multicolumn{1}{l}{0} & \multicolumn{1}{l}{0} \\
          & 60    & ABC   & \multicolumn{1}{l}{2.18e - 04} & \multicolumn{1}{l}{1.24e - 03} & \multicolumn{1}{l}{3.25e - 04} & \multicolumn{1}{l}{5.78e - 04} & \multicolumn{1}{l}{3.51e - 04} \\
          &       & MABC  & \multicolumn{1}{l}{2.55e - 16} & \multicolumn{1}{l}{1.90e - 15} & \multicolumn{1}{l}{4.45e - 16} & \multicolumn{1}{l}{8.20e - 16} & \multicolumn{1}{l}{4.69e - 16} \\
          &       & MPPA  & \multicolumn{1}{l}{0} & \multicolumn{1}{l}{0} & \multicolumn{1}{l}{0} & \multicolumn{1}{l}{0} & \multicolumn{1}{l}{0} \\
          & 100   & ABC   & \multicolumn{1}{l}{7.13e - 04} & \multicolumn{1}{l}{1.27e - 02} & \multicolumn{1}{l}{7.13e - 04} & \multicolumn{1}{l}{7.62e - 03} & \multicolumn{1}{l}{5.10e - 03} \\
          &       & MABC  & \multicolumn{1}{l}{2.38e - 15} & \multicolumn{1}{l}{9.68e - 15} & \multicolumn{1}{l}{7.76e - 15} & \multicolumn{1}{l}{5.83e - 15} & \multicolumn{1}{l}{1.97e - 15} \\
          &       & MPPA  & \multicolumn{1}{l}{0} & \multicolumn{1}{l}{0} & \multicolumn{1}{l}{0} & \multicolumn{1}{l}{0} & \multicolumn{1}{l}{0} \\
    17    & 30    & ABC   & \multicolumn{1}{l}{1.38e - 01} & \multicolumn{1}{l}{1.62e - 01} & \multicolumn{1}{l}{1.39e - 01} & \multicolumn{1}{l}{1.46e - 01} & \multicolumn{1}{l}{1.09e - 02} \\
          &       & MABC  & \multicolumn{1}{l}{0} & \multicolumn{1}{l}{0} & \multicolumn{1}{l}{0} & \multicolumn{1}{l}{0} & \multicolumn{1}{l}{0} \\
          &       & MPPA  & \multicolumn{1}{l}{0} & \multicolumn{1}{l}{0} & \multicolumn{1}{l}{0} & \multicolumn{1}{l}{0} & \multicolumn{1}{l}{0} \\
          & 60    & ABC   & \multicolumn{1}{l}{1.87e - 01} & \multicolumn{1}{l}{3.51e - 01} & \multicolumn{1}{l}{2.93e - 01} & \multicolumn{1}{l}{2.77e - 01} & \multicolumn{1}{l}{6.77e - 02} \\
          &       & MABC  & \multicolumn{1}{l}{0} & \multicolumn{1}{l}{1.42e - 14} & \multicolumn{1}{l}{7.10e - 15} & \multicolumn{1}{l}{9.94e - 15} & \multicolumn{1}{l}{5.68e - 15} \\
          &       & MPPA  & \multicolumn{1}{l}{0} & \multicolumn{1}{l}{0} & \multicolumn{1}{l}{0} & \multicolumn{1}{l}{0} & \multicolumn{1}{l}{0} \\
          & 100   & ABC   & \multicolumn{1}{l}{6.66e - 00} & \multicolumn{1}{l}{7.48e - 00} & \multicolumn{1}{l}{7.26e - 00} & \multicolumn{1}{l}{7.07e - 00} & \multicolumn{1}{l}{4.08e - 00} \\
          &       & MABC  & \multicolumn{1}{l}{4.26e - 14} & \multicolumn{1}{l}{5.68e - 14} & \multicolumn{1}{l}{4.26e - 14} & \multicolumn{1}{l}{5.21e - 14} & \multicolumn{1}{l}{6.69e - 15} \\
          &       & MPPA  & \multicolumn{1}{l}{0} & \multicolumn{1}{l}{0} & \multicolumn{1}{l}{0} & \multicolumn{1}{l}{0} & \multicolumn{1}{l}{0} \\
    18    & 30    & ABC   & \multicolumn{1}{l}{4.147e - 01} & \multicolumn{1}{l}{4.598e - 01} & \multicolumn{1}{l}{4.524e - 01} & \multicolumn{1}{l}{4.413e - 01} & \multicolumn{1}{l}{1.81e - 02} \\
          &       & MABC  & \multicolumn{1}{l}{2.277e - 01} & \multicolumn{1}{l}{3.455e - 01} & \multicolumn{1}{l}{3.121e - 01} & \multicolumn{1}{l}{2.952e - 01} & \multicolumn{1}{l}{3.17e - 02} \\
          &       & MPPA  & \multicolumn{1}{l}{3.14e - 03} & \multicolumn{1}{l}{9.57e - 02} & \multicolumn{1}{l}{3.36e - 02} & \multicolumn{1}{l}{3.81e - 02} & \multicolumn{1}{l}{2.65e - 02} \\
          & 60    & ABC   & \multicolumn{1}{l}{4.960e - 01} & \multicolumn{1}{l}{4.976e - 01} & \multicolumn{1}{l}{4.974e - 01} & \multicolumn{1}{l}{4.971e - 01} & \multicolumn{1}{l}{5.90e - 04} \\
          &       & MABC  & \multicolumn{1}{l}{4.796e - 01} & \multicolumn{1}{l}{4.903e - 01} & \multicolumn{1}{l}{4.850e - 01} & \multicolumn{1}{l}{4.840e - 01} & \multicolumn{1}{l}{3.62e - 03} \\
          &       & MPPA  & \multicolumn{1}{l}{9.29e - 03} & \multicolumn{1}{l}{1.98e - 01} & \multicolumn{1}{l}{5.13e - 02} & \multicolumn{1}{l}{6.84e - 02} & \multicolumn{1}{l}{4.71e - 02} \\
          & 100   & ABC   & \multicolumn{1}{l}{4.996e - 01} & \multicolumn{1}{l}{4.998e - 01} & \multicolumn{1}{l}{4.998e - 01} & \multicolumn{1}{l}{4.997e - 01} & \multicolumn{1}{l}{4.51e - 05} \\
          &       & MABC  & \multicolumn{1}{l}{4.988e - 01} & \multicolumn{1}{l}{4.992e - 01} & \multicolumn{1}{l}{4.991e - 01} & \multicolumn{1}{l}{4.990e - 01} & \multicolumn{1}{l}{1.75e - 04} \\
          &       & MPPA  & \multicolumn{1}{l}{0} & \multicolumn{1}{l}{1.77e - 01} & \multicolumn{1}{l}{7.85e - 02} & \multicolumn{1}{l}{8.79e - 02} & \multicolumn{1}{l}{5.33e - 02} \\
    \bottomrule
    \end{tabular}%
  \label{tab:addlabel}%
\end{table*}%
}
\subsection{Conclusion}

Optimisation problems are becoming more and more complex and unavoidable in most human activities. Although a variety of algorithms and heuristics to deal with them have been developed, new approaches are needed as the size and complexity of these problems increase. In recent years, heuristics and in particular those inspired by Nature, are becoming more efficient and robust. We have designed a Nature-inspired algorithm based on the way plants and in particular the strawberry plant, propagate. The original PPA algorithm has been only summarily tested on low dimension problems to establish its credentials. Surprisingly, despite being very simple and requiring few parameters, it managed to solve those problems rather well, \cite{Salhi2010PPA}. In this paper, we have presented a modified version of PPA, which is referred to as MPPA. The improvements concern the way new solutions at the end of runners are calculated, i.e. Equations \ref{eq:perturbation1}-\ref{eq:perturbation3}. The resulting algorithm has been tested on a more extensive test bench with a large number of functions having interesting characteristics such as multi-modality and non-separability in high dimensions, up to a 100. The results show that MPPA outperforms ABC and its more robust modification MABC on most of the test functions. In conclusion, MPPA provides us with a robust, easy to implement method for nonlinear, non-convex high dimensional optimisation problems.

\section*{Acknowledgments}

This work is supported by Abdul Wali Khan University, Mardan, Pakistan, Grant No. F.16-5/ P\& D/ AWKUM /238.


\newpage
\bibliographystyle{plain}
\bibliography{Sulaimanbib}

\end{document}